\documentclass[12pt]{article}
\usepackage{amsfonts}
\usepackage{amsmath}
\usepackage{amssymb}
\usepackage{amscd}
\usepackage{eepic}
\usepackage{epic}
\usepackage[francais]{babel}
\usepackage[latin1]{inputenc}
\usepackage{color}
\usepackage{graphics}
\begin{document}
\newtheorem{Theoreme}{Th\'eor\`eme}[section]
\newtheorem{Theorem}{Theorem}[section]
\newtheorem{Th}{Th\'eor\`eme}[section]
\newtheorem{De}[Th]{D\'efinition}
\newtheorem{Pro}[Th]{Proposition}
\newtheorem{Lemma}[Theorem]{Lemma}
\newtheorem{Proposition}[Theoreme]{Proposition}
\newtheorem{Lemme}[Theoreme]{Lemme}
\newtheorem{Corollaire}[Theoreme]{Corollaire}
\newtheorem{Consequence}[Theoreme]{Cons\'equence}
\newtheorem{Remarque1}[Theoreme]{Remarque}
\newtheorem{Convention}[Theoreme]{{\sc Convention}}
\newtheorem{PP}[Theoreme]{Propri\'et\'es}
\newtheorem{Conclusion}[Theoreme]{Conclusion}
\newtheorem{Ex}[Theoreme]{Exemple}
\newtheorem{Definition}[Theoreme]{D\'efinition}
\newtheorem{Remark1}[Theorem]{Remark}
\newtheorem{Not}[Theoreme]{Notation}
\newtheorem{Nota}[Theorem]{Notation}
\newtheorem{Propo}[Theorem]{Proposition}
\newtheorem{exercice1}[Th]{Lemme-Confié au lecteur}
\newtheorem{Corollary}[Theorem]{Corollary}
\newtheorem{PPtes}[Th]{Propri\'et\'es}
\newtheorem{Defi}[Theorem]{Definition}
\newtheorem{Example1}[Theorem]{Example}
\newenvironment{Proof}{\medbreak{\noindent\bf Proof }}{~{\hfill $\bullet$\bigbreak}}

\newenvironment{Demonstration}{\medbreak{\noindent\bf D\'emonstration
 }}{~{\hskip 3pt$\bullet$\bigbreak}} 

\newenvironment{Remarque}{\begin{Remarque1}\em}{\end{Remarque1}} 
\newenvironment{Remark}{\begin{Remark1}\em}{\end{Remark1}}
\newenvironment{Exemple}{\begin{Ex}\em}{~{\hskip
3pt$\bullet$}\end{Ex}} 
\newenvironment{exercice}{\begin{exercice1}\em}{\end{exercice1}}
\newenvironment{Notation}{\begin{Not}\em}{\end{Not}}
\newenvironment{Notation1}{\begin{Nota}\em}{\end{Nota}}
\newenvironment{Example}{\begin{Example1}\em}{~{\hskip
3pt$\bullet$}\end{Example1}} 

\newenvironment{Remarques}{\begin{Remarque1}\em \ \\* }{\end{Remarque1}}

\renewcommand{\Re}{{\cal R}}
\newcommand{\Dim}{{\rm Dim\,}}
\renewcommand{\Im}{{\cal F}}
\newcommand{\finpreuve}{~{\hskip 3pt$\bullet$\bigbreak}}
\newcommand{\hp}{\hskip 3pt}
\newcommand{\hph}{\hskip 8pt}
\newcommand{\hphh}{\hskip 15pt}
\newcommand{\vp}{\vskip 3pt}
\newcommand{\vpv}{\vskip 15pt}
\newcommand{\IP}{{\mathbb{IP}}}
\newcommand{\rd}{{\mathbb R}^2}
\newcommand{\R}{{\mathbb R}}

\newcommand{\Hyper}{{\mathbb H}}
\newcommand{\Int}{{\mathbb I}}
\newcommand{\Boule}{{\mathbb B}(0,1)}
\newcommand{\Cantor}{{\mathbb K}}
\newcommand{\Sp}{{\mathbb S}}
\newcommand{\K}{{\mathbb K}}
\newcommand{\B}{{\mathbb B}}
\newcommand{\Ho}{{\mathbb H}}
\newcommand{\Nat}{{\mathbb N}}
\newcommand{\N}{{\mathbb N}}
\newcommand{\Proba}{{\mathbb P}}
\newcommand{\Esp}{{\mathbb E}}
\newcommand{\Complex}{{\mathbb C}}
\newcommand{\Ha}{{\cal H}}
\newcommand{\Harm}{{\bold H}}
\newcommand{\Lcal}{{\cal L}}
\newcommand{\ds}{\displaystyle}
\newcommand{\un}{\bold 1}
\newcommand{\Cone}{C(x,r,\epsilon ,\Phi)}
\newcommand{\Cn}{C(x,2^{-n},\epsilon ,\Phi )}
\newcommand{\Tranche}{W(x,r,\epsilon,\Phi)}
\newcommand{\Wn}{W(x,2^{-n},\epsilon,\Phi)}
\newcommand{\WFn}{W(x,2^{-n},\epsilon,\Phi)\cap F}
\newcommand{\ovec}{\overrightarrow}
\newcommand{\red}{{\bold R}}
\newcommand{\dimH}{\dim_{\Ha}}
\newcommand{\diam}{\mbox{diam}}
\newcommand{\diamit}{\mbox{\em diam}}
\newcommand{\para}{\vskip 2mm}
\newcommand{\cod}{\stackrel{\mbox{\tiny cod}}{\sim}}
\newcommand{\cardit}{\mbox{\em card}}
\newcommand{\card}{\mbox{card}}
\newcommand{\Sphere}{{\mathbb S}_d}
\newcommand{\distit}{\mbox{\em dist}}
\newcommand{\Tri}{{\cal P}}
\newcommand{\LL}{{\mathcal L}}
\newcommand{\infess}{\mbox{inf\,ess}}
\newcommand{\supess}{\mbox{sup\,ess}}

\definecolor{darkblue}{rgb}{0,0,.5}
\def\u{\underline }
\def\o{\overline}
\def\h{\hskip 3pt}
\def\hh{\hskip 8pt}
\def\hhh{\hskip 15pt}
\def\v{\vskip 8pt}
\def\vv{\vskip 15pt}
\font\courrier=cmr12
\font\grand=cmbxti10
\font\large=cmbx12
\font\largeplus=cmr17
\font\small=cmbx8
\font\nor=cmbxti10
\font\smaller=cmr8
\font\smallo=cmbxti10
\openup 0.3mm
%
%
%

\title{\color{blue} Multifractal Analysis of inhomogeneous Bernoulli products}
\author{Athanasios BATAKIS and Benoît TESTUD}
\date{}
\maketitle
{\bf Abstract} We are interested to the multifractal analysis of  inhomogeneous Bernoulli products  which are also known as coin tossing measures. We give conditions ensuring the validity of the  multifractal formalism for such measures. On another hand, we show that these measures can have a dense set of phase transitions. 
 

\section{Introduction}

Let us consider the dyadic tree $\mathbb T$ (even though all the results in this paper can be easily generalised to any $\ell$-adic structure, $\ell\in\N$), let $\Sigma=\{0,1\}^{\N}$ be its limit (Cantor) set and denote by $(\Im_n)_{n\in\N}$ the associated filtration with the usual $0-1$ encoding. 

For $\epsilon_1,...,\epsilon_n\in\{0,1\}$ we denote by $I_{\epsilon_1...\epsilon_n}$ the cylinder of the $n$th generation defined by 
$I_{\epsilon_1...\epsilon_n}=\{x=(i_1,...,i_n,i_{n+1},...)\in \Sigma,\; ;\;i_1=\epsilon_1,...,i_n=\epsilon_n\}$. For every $x\in \Sigma$, $I_n(x)$ stands for   the cylinder  of  $\mathcal{F}_n$ containing~$x$.

If $(p_n)_n$ is a sequence of weights, $p_i\in (0,1)$,  we are interested in Borel measures $\mu$ on $\sigma$  defined 
in the following way 
\begin{equation}\label{inv}
\mu(I_{\epsilon_1...\epsilon_n})=\prod_{j=1}^n p_j^{1-\epsilon_j}(1-p_j)^{\epsilon_j}.
\end{equation}

This type of measure will be referred to as an {\em inhomogeneous Bernoulli product}.
The aim of this paper is to study multifractal properties of such measures.

The particular case where the sequence $(p_n)$ is constant is well-known and provides an example of measure satisfying the multifractal formalism (see e.g \cite{fal}). In the general case,    Bisbas   in \cite{Bis} gave a sufficient condition on the sequence $(p_n)$  ensuring that $\mu$ is a multifractal measure (i.e. the level sets are not empty) . However, the work of Bisbas does not provide the dimension of the level sets $E_\alpha$ associated to  the measure $\mu$. 

Let us give a brief description of multifractal formalism. 
For  a probability measure $m$ on $\Sigma$,  we  define the {\em local dimension} (also called   H\"older exponent) of $m$ at $x\in \Sigma$ by $$\displaystyle \alpha(x)=\liminf_{n\rightarrow +\infty} \alpha_n(x)=\liminf_{n\rightarrow +\infty} -\frac{\log m(I_n(x))}{n\log 2}.$$   The  aim of  multifractal analysis is  to find the Hausdorff dimension, $\dim(E_\alpha)$,  of  the  level  set  $ E_\alpha=\left\{x : \alpha(x)=\alpha\right\}$ for  $\alpha>0$.   The function  $f(\alpha)=\dim(E_\alpha)$ is called the {\em singularity spectrum} (or multifractal spectrum) of $m$ and we say that $m$ is a {\em multifractal measure} when $f(\alpha)>0$ for several $\alpha{}'s$.

The concepts underlying the  multifractal decomposition of a measure  go back to an early  paper of    Mandelbrot \cite{mandel}.   In the 80's multifractal measures were used    by physicists to study various  models arising from natural phenomena.   In    fully developped turbulence they  were  used by Frisch and Parisi \cite{fp} to investigate the intermittent behaviour in the regions of high vorticity. In dynamical system theory they were used by Benzi et al.~\cite{ben}  to measure how often a given region of the attractor is visited.  In diffusion-limited aggregation (DLA) they were  used by Meakin et al.~\cite{mea}  to describe the probability of a random walk landing to the neighborhood of a given site on the aggregate. 

In order to determine the function $f(\alpha)$,       Hentschel and  Procaccia  \cite{hp}  used  ideas based on  Renyi entropies \cite{re} to  introduce the generalized dimensions  $D_q$ defined by  \[D_q= \lim_{n\rightarrow +\infty} \frac{1}{q-1}  \frac{\log\left(\sum_{I\in\mathcal{F}_n} m(I)^q\right)}{n \log 2},\] (see also \cite{grp,gr}). From a physical and heuristical point of view, Halsey et al. \cite{has} showed that the singularity spectrum $f(\alpha)$ and the generalized dimensions $D_q$ can be derived from each other. The Legendre transform  turned out to be a useful tool linking  $f(\alpha)$ and $D_q$.  More precisely,  it was suggested that   
\begin{eqnarray}\label{fm} f(\alpha)=\dim(E_{\alpha})=\tau^*(\alpha)=\inf(\alpha q+\tau(q),\ q\in \R),\end{eqnarray}  where \[\tau (q) = \limsup_{n\rightarrow +\infty}\,\tau _n (q)\quad\mbox{ with }\quad \tau _n (q)=\frac{1}{n\log 2}\log\left(\sum_{I\in\mathcal{F}_n} m(I)^q\right).\]   (The sum  runs over the cylinders  $I$  such that $m(I)\not=0$.)   The  function $\tau(q)$ is  called the  $L^q$-spectrum of $m$   and if the limit exists $\tau(q)=(q-1)D_q$.

Relation \eqref{fm} is called the multifractal formalism and in many aspects it is analogous to the well-known thermodynamic formalism developed by Bowen \cite{bo} and Ruelle \cite{ru}. 

For number of measures, relation \eqref{fm}    can  be verified rigorously. In particular, if the sequence $(p_n)$ is constant or periodic, the measure $\mu$ given by  \eqref{inv}  satisfies  the multifractal formalism (e.g. \cite{fal}). Moreover some rigourous results have already been obtained for some invariant measures in dynamical systems (e.g \cite{col, Fan,rand}), for some self-similars measures under separation conditions (e.g \cite{cm,lng,ol})  and for quasiindependent measures(e.g \cite{BMP,heurt,test}).

     The  minoration of $\dim(E_\alpha)$ usually follows on the
existence of a shift-invariant and ergodic measure     $m_q$ (the so-called    {\em Gibbs measure} \cite{mi}),   satisfying  
\begin{eqnarray*}\forall n,\,\,\forall I\in \mathcal{F}_n, \quad
  \frac1C  m(I)^q 2^{-n\tau(q)}\le  m_q(I)\le C m(I)^q
  2^{-n\tau(q)},\end{eqnarray*} where the constant  $C>0$ is
independent of $n$ and $I$. If  $\tau$ is differentiable at $q$,   the measure  $m_q$  is 
supported by   $E_{-\tau'(q)}$ and      Brown, Michon
and Peyri\`ere    established  \cite{BMP,p}  that
\begin{eqnarray*}
\dim(E_{-\tau'(q)})
=\tau^*(-\tau'(q))=-q\tau'(q)+\tau(q).\end{eqnarray*}
 
If  the weights $p_n$ are not all  the same,  the measure $\mu$ is in general no shift-invariant and we cannot apply  classical tools  of ergodic theory, as Shannon-McMillan theorem (e.g \cite{Bil}), to get a lower bound of $\dim(E_\alpha)$.

 Let us introduce the other following level sets defined by
$$\underline E_{\alpha}=\left\{x\; ;\;
\alpha(x)\le \alpha\right\},  \;   \overline F_{\alpha}=\left\{x\; ;\;
\limsup_{n\to\infty}\alpha_n(x)\ge \alpha\right\},$$ and $$F_{\alpha}=\left\{x\; ;\;
\limsup_{n\to\infty}\alpha_n(x)= \alpha\right\}.$$
 
We can now state our main results. In section 2, we prove the following.  
\begin{Theorem}\label{1}
Let $\mu$ be an inhomogeneous Bernoulli product on $\Sigma$ and $q\in\R$.
 We have
$$
\liminf_{n\to\infty}-q\tau_{\mu,n}'(q)+\tau_{\mu,n}(q)\le\dim\left(\underline E_{-\tau'(q^-)}\cap \overline F_{-\tau'(q^+)}\right)\le
\inf\left\{\tau^*(-\tau'(q^+)), \tau^*(-\tau'(q^-))\right\}.
$$  \end{Theorem}

The proof of the lower bound relies on the construction of a special  inhomogeneous  Bernoulli product which has the dimension of the level set studied.

In section 3 we are interested to the case where the sequence $\tau_{\mu,n}(q)$ converges. In this situation, we prove that the multifractal formalism holds for  $\alpha=-\tau_\mu'(q)$--if it exists. More precisely, we have   
\begin{Theorem}\label{4}
Suppose that the sequence $(\tau_{\mu, n}(q))$ converges at a point $q>0$. If $\tau_{\mu}'(q)$ exists and if $\alpha=-\tau_\mu'(q)$, we have 
\begin{equation}
\dim \left(E_{\alpha}\cap F_{\alpha}\right)=\tau_\mu^*(\alpha)=\alpha q+\tau_\mu(q).
\end{equation}
\end{Theorem}

Theorem \ref{4} lead us to study the differentiability  of the $L^q$-spectrum $\tau_\mu(q)$. A point $q$ will be called a {\em phase transition} if $\tau_\mu'(q)$ does not exist.  In  section 4, we are interested to the existence of   phase transitions.  More precisely, we prove the following.  
\begin{Theorem}\label{3}
There exist inhomogeneous Bernoulli products $\mu$  presenting a dense set of 
phase transitions.
\end{Theorem}
\section{Proof of theorem \ref{1}}
We begin by a preliminary result.
\begin{Lemma}\label{secderbdd}
If $\mu$ is an inhomogeneous Bernoulli product, then $(\tau_{\mu,n}'')$
 are locally uniformly bounded on $(0,+\infty)$.
 \end{Lemma} 
 \begin{Proof}
 We denote by $\beta(p_i)$ the Bernoulli homogeneous measure of parameter $p_i$ and by $\tau(p_i,q)$ it's $\tau$ function, 
 $\tau(p_i,q)=\log(p_i^q+(1-p_i)^q)$. Using the fact that $\mu$ is the product of $\beta(p_i)$ we easily obtain
 $$\tau_{\mu,n}(q)=\frac1n\sum_{i=0}^n\tau(p_i,q)\;\,\;q>0.$$ 
 It is therefore sufficient to show that ,for any $q_0>0$, there exists a constant $C=C(q_0)$ such that 
for all $p\in(0,1)$ and   all $q>q_0$, $\displaystyle\frac{\partial^2\tau(p,q)}{\partial q^2}\le C$.
 The proof is straigthforward:
\begin{eqnarray*}
\frac{\partial^2\tau(p,q)}{\partial q^2}&=&
\frac{\left(p^q(\log p)^2+(1-p)^q(\log(1-p))^2\right)}{(p^q+(1-p)^q)}-\frac{\left(p^q\log p+(1-p)^q\log(1-p)\right)^2}{(p^q+(1-p)^q)^2}\\
&=&\frac{p^q(1-p)^q\left((\log p)^2+(\log(1-p))^2-2\log p\log(1-p)\right)}{(p^q+(1-p)^q)^2}\\
&=&\frac{p^q(1-p)^q\left(\log\frac{p}{1-p}\right)^2}{(p^q+(1-p)^q)^2}\le [4p(1-p)]^q(\log p)^2\le[4p(1-p)]^{q_0}(\log p)^2,
\end{eqnarray*}
which is uniformly bounded on $p\in(0,1)$ and the proof is complete.
\end{Proof}

Lemma \ref{secderbdd} allows us to give estimates for the lower and the upper Hausdorff dimension of the measure $\mu$. They are respectively defined by 
$$\dim_*(\mu)=\inf\{\dim(E),\,\,\, \mu(E)>0\} \,;\,\,\dim^*(\mu)=\inf\{\dim(E),\,\,\, \mu(E)=1\}.$$  We say that $\mu$ is exact if $\dim_*(\mu)=\dim^(\mu)$ and we note $\dim(\mu)$ the common value. In the same way,  we can define  the lower and the upper packing  dimension of the measure $\mu$.  It is well known that there exist some relations between these quantities and the derivatives of the function $\tau_\mu(q)$ at $q=1$. More precisely, it is proved in \cite{Fan,heurt} that 
$$-\tau_\mu'(1+)\le\dim_*(\mu)\le  h_*(\mu)\le h^*(\mu)\le \Dim^*(\mu)\le -\tau_\mu'(1-),$$ where $ h_*(\mu)$ and $ h^*(\mu)$ stand for the lower and the upper entropy of the measure $\mu$, defined as 
$$ h_*(\mu)=\liminf -\frac{1}{n\log 2} \sum_{I\in\Im_n} \mu(I)\log \mu(I)= \liminf -\tau_{\mu_n}'(1)$$ and  $$ h^*(\mu)=\limsup -\frac{1}{n\log 2} \sum_{I\in\Im_n} \mu(I)\log \mu(I)= \limsup - \tau_{\mu_n}'(1).$$  

By Lemma \ref{secderbdd}, we deduce (see \cite{BH,heurt}) the following remark.
 \begin{Remark}\label{Heurt}
 If $\mu$ is an inhomogeneous Bernoulli product then 
$$\displaystyle\dim\mu=\liminf_{n\to\infty}-\tau_{\mu_n}'(1)=-\tau_{\mu}'(1^+)=h_*(\mu).$$ and 
$$\displaystyle\Dim\mu=\limsup_{n\to\infty}-\tau_{\mu_n}'(1)=-\tau_{\mu}'(1^-)=h^*(\mu).$$
 \end{Remark}

Fix $q\in \R$. To prove theorem \ref{1}, we construct an auxiliary  measure $\nu$ supported by the set $\underline E_{-\tau'(q^-)}\cap \overline F_{-\tau'(q^+)}$. More precisely, we consider a sequence of measures  $\nu_n$ satisfying 
$$\nu_n(I)=\frac{\mu(I)^q}{\sum_{I\in\Im_n}\mu(I)^q}=\mu(I)^q|I|^{\tau_{\mu,n}(q)},$$  if  $I\in\Im_n$. 
The following lemma implies that  the sequence $(\nu_n)$ converges in the weak$^*$ sense to a probability  measure  $\nu$ which 
is also an inhomogeneous Bernoulli product. 
 \begin{Lemma}\label{nuthing}
Let $n\in\N$ and $I\in \Im_n$. If  $\mu$ is an inhomogeneous Bernoulli product, we have 
$\nu_n(I) =\nu_{n+1}(I)$. \end{Lemma}
\begin{Proof}
Take $n>0$ and $I\in\Im_n$. We can compute
$$\nu_{n+1}(I)=\frac{\sum_{J\in\Im_1}\mu(IJ)^q}{\sum_{I\in\Im_n}\sum_{J \in\Im_1}\mu(IJ)^q}=
\frac{\mu(I)^q(p_{n+1}^q+(1-p_{n+1})^q)}{\sum_{I\in\Im_n}(p_{n+1}^q+(1-p_{n+1})^q)\mu(I)^q}$$
and therefore
$\nu_{n+1}(I)=\nu_n(I)$ for all $I\in \Im_n$.
\end{Proof}
By remark \ref{Heurt}, we then deduce that the Hausdorff and the packing dimension of $\nu$ are given by an entropy formula. In other terms, we have $$\displaystyle\dim\nu=\liminf_{n\to\infty}-\tau_{\nu,n}'(1)=h_*(\nu)$$
and
$$\displaystyle\Dim\nu=\limsup_{n\to\infty}-\tau_{\nu,n}'(1)=h^*(\nu).$$

Now we cam prove Theorem \ref{1}.

 \begin{Proof}{\bf of Theorem \ref{1}}
The upper bound is a well known fact of multifractal formalism (see for instance \cite{BMP}). 
In fact we have 
\begin{enumerate}
\item If $\alpha\le -\tau'(0^+)$ then 
$\dim E_{\alpha}\le\dim \underline E_{\alpha}\le\tau^*(\alpha).$
\item If $\alpha\ge -\tau'(0^-)$ then 
$\dim F_{\alpha}\le\dim \overline F_{\alpha}\le\tau^*(\alpha).$
\item $-\tau'(0^+)\le \alpha\le -\tau'(0^-)$ then $\tau^*(\alpha)=\tau(0)$ and the upper bound is trivial.
\end{enumerate}

Lemma \ref{nuthing} and a straightforward computation imply
$\tau_{\nu,n}(s)=\tau_{\mu,n}(qs)-s\tau_{\mu,n}(q)$.
Using once again the (inhomogeneous) Bernoulli property of $\mu$ and remark \ref{Heurt} we deduce that 
$$-\tau_{\nu}'(1^+)=\liminf-\tau_{\nu,n}'(1)=\liminf\left(-q\tau_{\mu,n}'(q)+\tau_{\mu,n}(q)\right).$$

The following lemma then implies the lower bound.
\begin{Lemma} We have $\nu\left(\underline E_{-\tau'(q^-)}\cap\overline F_{-\tau'(q^+)}\right)=1$.
\end{Lemma}
\begin{Proof}
For $\eta>0$ we put $\beta=-\tau_{\mu}'(q^-)+\eta$ and we prove that $\nu(\Sigma\setminus\underline E_{\beta})=0$ ; 
it can be shown in a similar way that
$\nu(\Sigma\setminus\overline F_{\gamma})=0$ for $\gamma<-\tau_{\mu}'(q^+)$. The lemma then easily follows.

It suffices to show that 
$\displaystyle \Sigma\setminus E_{\beta}=\left\{x\in\Sigma\; ;\; 
\liminf_{k\to\infty}\alpha_{n}(x)> \beta\right\}$ is of 0 $\nu$-measure. 
Consider the collection $\Re_{n}(\beta)$ of cylinders $I\in\Im_{n}$  satisfying
$\ds\frac{\log\mu(I)}{\log|I|}>\beta$. It is clear that $\ds \Sigma\setminus E_{\beta}=\limsup_{n\to\infty}\Re_{n}(\beta)$.

Let $(\tau_{\mu,n_k})_{k\in\N}$ be the subsequence of $(\tau_{\mu,n})_{n\in\N}$ such that $\lim_{k\to\infty}\tau_{\mu,n_k}(q)=\tau_{\mu}(q)$. 
Using the convergence of $\tau_{\mu,n_k}(q)$ we can choose
(and fix) $t<0$ such that for  $k$ big enough
$$\tau_{\mu}(q+t)-\tau_{\mu,n_k}(q)<-\left(\beta-\frac{\eta}{2}\right)t=\left(\tau_{\mu}'(q^-)-\frac{\eta}{2}\right) t$$
We get $\ds \mu(I)^{-t}|I|^{\beta t}\le 1$ and hence
\begin{eqnarray*}
\sum_{I\in\Re_{n_k}(\beta)}\nu(I)&=&\sum_{I\in\Re_{n_k}(\beta)}\mu(I)^q|I|^{\tau_{\mu,n_k}(q)}=
\sum_{I\in\Re_{n_k}(\beta)}\mu(I)^{q+t}|I|^{\tau_{\mu,n_k}(q)-\beta t}\mu(I)^{-t}|I|^{\beta t}\\
&\le&\sum_{I\in\Re_{n_k}(\beta)}\mu(I)^{q+t}|I|^{\tau_{\mu,n_k}(q)-\beta t}\le 
\sum_{I\in\Im_{n_k}}\mu(I)^{q+t}|I|^{\tau_{\mu}(q+t)-\frac{\eta}{2}t}\\
&\le& \sum_{I\in\Im_{n_k}}\mu(I)^{q+t}|I|^{\tau_{\mu,n_k}(q+t)}=1,
\end{eqnarray*}
where for the last inequality we used the fact that $\ds\tau_{\mu}(q+t)=\limsup\tau_{\mu,n}(q+t)$. 
It easily follows that $\ds\limsup_{k\to\infty}\sum_{I\in\Re_{n_k}(\beta)}\nu(I)=0$ and the lemma is proved.
\end{Proof}
The proof of Theorem \ref{1} is now completed.\end{Proof}

Let $f$ and $g$ be the functions defined by $f(t)=\dim \underline E_t$ and $g(t)=\dim \overline F_t$ . 
Obviously, $f$ is increasing and $g$ is decreasing. 
Recall that $t$ is a non-stationary point of a monotone function $h$ if $h(s)\not=h(t)$ for all $s\not=t$.

Since $E_{\alpha}=\underline E_{\alpha}\setminus\bigcup_{\beta< \alpha} \underline E_{\beta}$,  we deduce from theorem \ref{1} the following.
\begin{Remark}
If $\alpha=-\tau'(q^-)$ for $q>0$  is a non-stationary point of $f$ there a sequence of $q_m\le q$ such that $\alpha_m=-\tau'(q_m^-)$  
are non-stationary points of $f$ converging to $\alpha$ and
$$\liminf_{n\to\infty}-q_m\tau_{\mu,n}'(q_m)+\tau_{\mu,n}(q_m)\le \dim E_{\alpha_m}=\dim \underline E_{\alpha_m} \le \tau^*(\alpha_m).$$
If $\alpha=-\tau'(q^+)$ for $q<0$  is a non-stationary point of $g$ there a sequence of $q_m\ge q$ such that $\alpha_m=-\tau'(q_m^+)$  
are non-stationary points of $g$ converging to $\alpha$ and
$$\liminf_{n\to\infty}-q_m\tau_{\mu,n}'(q_m)+\tau_{\mu,n}(q_m)\le \dim F_{\alpha_m}=\dim \overline F_{\alpha_m} \le \tau^*(\alpha_m).$$

We conjecture  that under the same conditions on $\alpha$ we should also have $\dim E_{\alpha}=\dim \underline E_{\alpha}$ 
($\dim F_{\alpha}=\dim \overline F_{\alpha}$ respectively).
\end{Remark}

\section{Some conditions  ensuring  the validity of  multifractal formalism}
In this section we prove Theorem \ref{4}. We will use  the following result.

\begin{Propo}\label{2}
For  $q>0$ consider $(\tau_{\mu,n_k})$ 
the subsequence of $(\tau_{\mu,n})$ such that 
$$\displaystyle \lim_{k\to\infty}\tau_{\mu,n_k}(q)=\limsup_{n\to\infty}\tau_{\mu,n}(q).$$ 
Then if $q$ is a differentiability point of $\tau_{\mu}$ we have 
$$\lim_{k\to\infty}\tau_{\mu,n_k}'=\tau_{\mu}'(q).$$
\end{Propo}

\begin{Proof} 
The proposition is a immediate consequence of the following lemmas. 
\begin{Lemma}\label{1.1}
Under the assumptions of proposition  \ref{2}
$$\tau_{\mu}'(q^+)\ge\limsup_{k\to\infty}\tau_{\mu,n_k}'(q),$$ 
where $\tau_{\mu}'(q^+)$ stands for the right hand dérivative ot $\tau_{\mu}$ at $q$.
\end{Lemma}
On the other hand, we get
\begin{Lemma}\label{1.2}
Under the assumptions of proposition \ref{2}
$$\tau_{\mu}'(q^-)\le\liminf_{k\to\infty}\tau_{\mu,n_k}'(q),$$
where $\tau_{\mu}'(q^-)$ stands for the left hand dérivative of $\tau_{\mu}$ at $q$.
\end{Lemma}

\begin{Proof}{\bf of lemma \ref{1.1}}. Take $\epsilon>0$ and $\tilde q>q$ satisfying 

\begin{eqnarray*}\label{restrictions1}
&&\left|\frac{\tau_{\mu}(\tilde q)-\tau_{\mu}(q)}{\tilde q-q}-\tau_{\mu}'(q^+)\right|<\epsilon/8\\
&&|\tilde q-q|\sup_{n\in\N}||\tau_{\mu,n}''||_{\infty}<\epsilon/8.
\end{eqnarray*}
 and consider $(\tilde n_k)$ such that 
 $\displaystyle \lim_{k\to\infty}\tau_{\mu,\tilde n_k}(\tilde q)=
\limsup_{n\to\infty}\tau_{\mu,n}(\tilde q)$.

We can chose $k$ big enough to have 
\begin{eqnarray*}\label{restrictions2}
\frac{|\tau_{\mu,n_k}(q)-\tau_{\mu}(q)|}{|\tilde q-q|}&<&\epsilon/8\\
\frac{|\tau_{\mu,\tilde n_k}(\tilde q)-\tau_{\mu}(\tilde q)|}{|\tilde q-q|}&<&\epsilon/8\\
\tau_{\mu,n_k}(\tilde q)&\le& \tau_{\mu,\tilde n_k}(\tilde q)+(\tilde q-q)\epsilon/8.
\end{eqnarray*}
We then obtain
\begin{eqnarray*}
\tau_{\mu}'(q^+) &\ge& \displaystyle \frac{\tau_{\mu}(\tilde q)-\tau_{\mu}(q)}{\tilde q-q}-\epsilon/8
\ge \frac{\tau_{\mu,\tilde n_k}(\tilde q)-\tau_{\mu,n_k}(q)}{\tilde q-q}-\epsilon/4\\
&\ge&\displaystyle \frac{\tau_{\mu,n_k}(\tilde q)-\tau_{\mu,n_k}(q)}{\tilde q-q}-3\epsilon/8
\ge \tau_{\mu,n_k}'(q)-|\tilde q-q|\sup_{n\in\N}||\tau_{\mu,n}''||_{\infty}-\epsilon/2
\\&\ge& \tau_{\mu,n_k}'(q)-\epsilon.
\end{eqnarray*}
and the proof is completed.
\end{Proof} Lemma \ref{1.2} is proven in a similar manner and together with lemma \ref{1.1} 
provide the proposition's proof.
\end{Proof}

We can now  prove Theorem \ref{4}.
 
\begin{Proof}{\bf of Theorem \ref{4}}.
Let $\nu$ be the Gibbs-measure defined in lemma \ref{nuthing}. Since 
$$\tau_{\nu,n}(s)=\tau_{\mu,n}(qs)-s\tau_{\mu,n}(q)$$ we get 
$$\tau_{\nu,n}'(1)=q\tau_{\mu,n}'(q)-\tau_{\mu,n}(q).$$
Using the convergence of $\tau_{\mu,n}(q)$ we deduce from Proposition \ref{2} that
$$\lim_{n\to\infty}\tau_{\nu,n}'(1)=\lim_{n\to\infty}\left(q\tau_{\mu,n}'(q)-\tau_{\mu,n}(q)\right)=q\tau_\mu'(q)-\tau_\mu(q).$$
Lemma \ref{nuthing} then implies that $\tau_{\nu}'(1)$ exists and 
$$\dim\nu=\Dim\nu=-\tau_{\nu}'(1)=-q\tau_\mu'(q)+\tau_\mu(q).$$

On the other hand, for $I\in\Im_n$, we have 
$$\frac{\log\nu(I)}{\log |I|}=q\frac{\log\mu(I)}{\log|I|}+\tau_{\mu,n}(q)$$
 Since 
 $$\lim_{n\to\infty}\frac{\log\nu(I_n(x))}{\log |I_n(x)|}=\dim\nu=\Dim\nu\; \,;\nu\mbox{-a.s.}$$
 we obtain that $\ds \lim_{n\to\infty}\frac{\log\mu(I_n(x))}{\log |I_n(x)|}=-\tau_\mu'(q)$ , $\nu$-a.s.
 We conclude that $$\dim\left(E_{\alpha}\cap F_{\alpha}\right)\ge \dim\nu=\tau_\mu^*(\alpha).$$
 The opposite inequality being always valid, the proof  is done.
\end{Proof}
\section{Phase transitions}
\begin{Theorem}\label{cocorico} Let $\tau$ be a convex combination of functions $\tau(p_i,.)$ where $0<p_i\le 1/2$, $i=1,...,n$. For any $1<q_1<q_2<\infty$ there exists another convex combination $\tilde \tau$ of functions $\tau(p_j',.)$ such that
\begin{itemize}
\item $\tilde\tau(q_i)=\tau(q_i)$ and $\tilde\tau'(q_i)\not=\tau'(q_i)$, $i=1,2$,
\item for $q\in(q_1,q_2)$, $\tilde\tau(q)>\tau(q)$,
\item else, for $q\notin[q_1,q_2]$, $\tilde\tau(q)<\tau(q)$.
\end{itemize}
\end{Theorem}
\section{Proof of Theorem \ref{cocorico}} 
In this section whenever we use the notation $p_i$ for a weight in $(0,1)$ we will also note $\tau_i=\tau(p_i,.)$.
\begin{Lemma}\label{basic}
Take $\tau=\lambda\tau(p_1,.)+(1-\lambda)\tau(p_2,.)$ with $0<p_1<p_2<1/2$ and $\lambda\in(0,1)$. For $p_0\in(0,1/2)$
one of the following occurs:
\begin{enumerate}
\item  either $\tau(q)\not= \tau(p_0,q)$, for all $q>1$,
\item or, there exists $q_0>1$ such that $\tau(q)>\tau(p_0,q)$ for $q<q_0$ and $\tau(q)< \tau(p_0,q)$ for $q>q_0$. The point $q_0$ is, then, the unique point of equality between these functions.
\end{enumerate}
\end{Lemma}
To prove this lemma we need the following subsidiary result.
\begin{Lemma}\label{subsidiary} Let $p_1<p_2<p_3$ take values in  $(0,1/2)$ and $\tau_1,\tau_2,\tau_3$ be the functions $\tau(p_1,.),\tau(p_2,.),\tau(p_3,.)$ respectively. Then 
$\ds \frac{\tau_1-\tau_2}{\tau_2-\tau_3}$ is decreasing on $(1,+\infty)$.
\end{Lemma}
Although the proof only uses elementary calculus, it is a litte bit ``tricky'' and canot be omitted.
\begin{Proof}{\bf of Lemma \ref{subsidiary}} Taking into account the trivial equality
$$\tau(p',q)-\tau(p'',q)=\int_{p''}^{p'}\frac{\partial\tau}{\partial p}(p,q)dp$$ we only need to show that if $p'<p''$ then\; \;
$\ds\frac{\partial\tau}{\partial p}(p',q):\frac{\partial\tau}{\partial p}(p'',q)$\; \;is decreasing on $q\in(1,\infty)$.
We get
\begin{eqnarray*}
\ds\frac{\partial\tau}{\partial p}(p',q):\frac{\partial\tau}{\partial p}(p'',q) &=& \frac{1-(-1+1/p')^{q-1}}{1+(-1+1/p')^{q}}:\frac{1-(-1+1/p'')^{q-1}}{1+(-1+1/p'')^{q}}\\
&=& p''\frac{1-{s_1}^{q-1}}{1+{s_1}^{q}}:p'\frac{1-{s_2}^{q-1}}{1+{s_2}^{q}}
\end{eqnarray*}
where $s_1=-1+1/p'>1$ and $s_2=-1+1/p''>1$.

If we set $\ds f(s,q)=\ln\frac{1-{s}^{q-1}}{1+{s}^{q}}$, with $s,q>1$, it is  sufficient to prove that $\ds\frac{\partial f}{\partial s}f(s,q)$ is decreasing in $q$.
We calculate
$$\ds\frac{\partial f}{\partial s}f(s,q)=\frac{(q-1)s^{q-2}}{s^{q-1}-1}-\frac{qs^{q-1}}{s^{q}+1}.$$
We multiply by $s$ and need to show that
$\ds \frac{(q-1)s^{q-1}}{s^{q-1}-1}-\frac{qs^{q}}{s^{q}+1}$ is decreasing  which is equivalent to 
$q-1+\frac{q-1}{s^{q-1}-1}-q+\frac{q}{s^{q}+1}$ being decreasing. 

Put $Q=q-1$ ; it remains to show that
$\ds \frac{q-1}{s^{q-1}-1}+\frac{q}{s^{q}+1}=\frac{Q}{s^Q-1}+\frac{Q}{s^{Q+1}+1}+\frac{1}{s^{Q+1}+1}$ decreases in $Q>0$ and since the last term is decreasing it suffices to show that $\ds\frac{Q}{s^Q-1}+\frac{Q}{s^{Q+1}+1}$ is doing the same.
By taking derivatives we need to show that 
$$(s^Q-1)(s^{Q+1}+1)-s^Q\ln s^Q(s^{Q+1}+1)-s^{Q+1}\ln s^Q(s^{Q}-1)$$ is negative for $Q>0$, which is trivial since $\ds s^Q\ln s^Q>s^{Q}-1$.
\end{Proof}
\begin{Proof}{\bf of lemma \ref{basic}.}
Let us first remark that  $\tau$ and $\tau(p_0,.)$ can coincide at one point  only if $p_0\in(p_1,p_2)$. Moreover,
$\tau(q)=\tau(p_0,q)$ implies
$$\frac{\tau(p_1,q)-\tau(p_0,q)}{\tau(p_0,q)-\tau(p_2,q)}=\frac{\lambda}{1-\lambda}.$$
By lemma \ref{subsidiary} this can only occur once. The lemma \ref{basic} easily follows on the decreasing property of the ratio.
\end{Proof}

The following two lemmas prove Theorem \ref{cocorico} in the particular case $n=2$. 
\begin{Lemma}\label{system} 
Take $\lambda_1,\lambda_2\in(0,1)$ such that $\lambda_1+\lambda_2=1$ , $1<p_1<p_2<1/2$ and 
set $\tau=\lambda_1\tau_1+\lambda_2\tau_2$. Fix $1<q_1<q_2<+\infty$ and consider $p_1<p_4<p_2<p_5<1/2$ such that $\tau(p_4,q)=\tau(q)$. Then there is a unique convex combination $\tilde \tau$ of $\tau_1,\tau_4$ and $\tau_5$ such that
$$\tilde\tau(q_1)=\tau(q_1)\mbox{ and }\tilde\tau(q_2)=\tau(q_2).$$
Furthermore, for $i=1,2$, we have $\tau'(q_i)\not=\tilde\tau'(q_i)$ and $\tau(q)\not=\tilde\tau(q)$ if $q\not=q_i$.
\end{Lemma}
\begin{Proof}
It suffices to  show that the linear system 
\begin{equation*}\left\{\begin{array}{llllcl}
&\lambda_3\tau_1(q_1)&+\lambda_4\tau_4(q_1)&+\lambda_5\tau_5(q_1)&=&\tau(q_1)\\
&\lambda_3\tau_1(q_2)&+\lambda_4\tau_4(q_2)&+\lambda_5\tau_5(q_2)&=&\tau(q_2)\\
&\lambda_3&+\lambda_4&+\lambda_5&=&1
\end{array} \right.\hspace{1cm}{ \bf (S)}\end{equation*}
has a unique positive solution $(\lambda_3,\lambda_4,\lambda_5)$.  The existence of a unique solution is easy to verify. Let us show that this solution is positive.

    Since $\tau(q_1)=\tau_4(q_1)$, we have 
$$\lambda_3\tau_1(q_1)+\lambda_5\tau_5(q_1)=(1-\lambda_4)(\lambda_1\tau_1(q_1)+\lambda_2\tau_2(q_1))$$
which is equivalent to  
\begin{eqnarray}\label{s}\frac{\lambda_3}{\lambda_3+\lambda_5}\tau_1(q_1)+\frac{\lambda_5}{\lambda_3+\lambda_5}\tau_5(q_1)=\lambda_1\tau_1(q_1)+\lambda_2\tau_2(q_1).\end{eqnarray}
This implies   that $\lambda_3 \lambda_5>0$.  Moreover,  since $\tau_5<\tau_2$, we also  have $\frac{\lambda_3}{\lambda_3+\lambda_5} >\lambda_1$. 

Let us show that $\lambda_3$ and $\lambda_5$ are positive. Otherwise,  by the above remark, we have $\lambda_3<0$, $\lambda_5<0$ and $\lambda_4> 0$. By the system (S) we have 
\begin{eqnarray*}\tau_4(q)=\frac{\lambda_1-\lambda_3}{\lambda_4}\tau_1(q)+\frac{\lambda_2}{\lambda_4}\tau_2(q)-\frac{\lambda_5}{\lambda_4}\tau_5(q)\end{eqnarray*} at the points $q=q_1$ and $q=q_2$. We then obtain that  $$\frac{\lambda_1-\lambda_3}{\lambda_4}\frac{\tau_1-\tau_4}{\tau_4-\tau_2}(q)=\frac{\lambda_2}{\lambda_4}+-\frac{\lambda_5}{\lambda_4}\frac{\tau_4-\tau_5}{\tau_4-\tau_2}(q)$$ for $q=q_1$ and $q=q_2$. Since $p_1<p_4<p_2$, by Lemma \ref{subsidiary} the function $\frac{\tau_1-\tau_4}{\tau_4-\tau_2}$ is decreasing. On the other hand,  since  $p_4<p_2<p_5$,  Lemma \ref{subsidiary}  implies that the function $\frac{\tau_4-\tau_5}{\tau_4-\tau_2}=1+\frac{\tau_2-\tau_5}{\tau_4-\tau_2}$ is increasing. Thus, these functions cannot coincide at two  points so   we conclude  that $\lambda_3$ and $\lambda_5$ are positive.

Let us now prove that $\lambda_4>0$. By \eqref{s} we have $$\frac{\lambda_3}{\lambda_3+\lambda_5}\tau_1(q_1)+\frac{\lambda_5}{\lambda_3+\lambda_5}\tau_5(q_1)=\lambda_1\tau_1(q_1)+\lambda_2\tau_2(q_1)$$ which gives that $$\lambda_2\tau_2(q_1)=\left(\frac{\lambda_3}{\lambda_3+\lambda_5}-\lambda_1\right)\tau_1(q_1)+\frac{\lambda_5}{\lambda_3+\lambda_5}\tau_5(q_1).$$ Using Lemma \ref{subsidiary}, for $q> q_1$   we get
$$\lambda_2\tau_2(q)> \left(\frac{\lambda_3}{\lambda_3+\lambda_5}-\lambda_1\right)\tau_1(q)+\frac{\lambda_5}{\lambda_3+\lambda_5}\tau_5(q)$$ and $$\lambda_1\tau_1(q)+\lambda_2\tau_2(q)> \frac{\lambda_3}{\lambda_3+\lambda_5}\tau_1(q)+\frac{\lambda_5}{\lambda_3+\lambda_5}\tau_5(q).$$ In particular, for $q=q_2$ we find that $$\lambda_3\tau_1(q_2)+\lambda_5\tau_5(q_2)+\lambda_4\tau(q_2)<\tau_(q_2)=\lambda_3\tau_1(q_2)+\lambda_4\tau_4(q_2)+\lambda_5\tau_5(q_2)$$ and we deduce that $$\lambda_4 \tau(q_2)<\lambda_4 \tau_4(q_2).$$ It follows from Lemma \ref{subsidiary} $\lambda_4>0$.

The last assertion follows directly from the independancy of the vector families
$$\left\{\left(\begin{array}{c}\tau_1(q_1)\\ \tau_4(q_1)\\ \tau_5(q_1)\end{array}\right), \left(\begin{array}{c}\tau_1(q_2)\\ \tau_4(q_2)\\ \tau_5(q_2)\end{array}\right), \left(\begin{array}{c}\tau_1'(q_i)\\ \tau_4'(q_i)\\ \tau_5'(q_i)\end{array}\right)\right\} $$ and
$$\left\{\left(\begin{array}{c}\tau_1(q_1)\\ \tau_4(q_1)\\ \tau_5(q_1)\end{array}\right), \left(\begin{array}{c}\tau_1(q_2)\\ \tau_4(q_2)\\ \tau_5(q_2)\end{array}\right), \left(\begin{array}{c}\tau_1(q)\\ \tau_4(q)\\ \tau_5(q)\end{array}\right)\right\} ,$$
which can be easily established.
\end{Proof}
\begin{Lemma}
The functions $\tau$ and $\tilde\tau$ defined in lemma \ref{system} verify
 $\tilde\tau(q)>\tau(q)$ if and only if $q\in(q_1,q_2)$.
\end{Lemma}

\begin{Proof}
Let us first remark that for $\lambda_3$, $\lambda_4$ and $\lambda_5$ defined by the linear system $({\bf S})$ we have
$\ds\frac{\lambda_3}{1-\lambda_4}\tau_1(q_1)+\frac{\lambda_5}{1-\lambda_4}\tau_5(q_1)=\tau(q_1)=\tau_4(q_1)$. Put $\rho= \frac{\lambda_3}{\lambda_3+\lambda_5}\tau_1+\frac{\lambda_5}{\lambda_3+\lambda_5}\tau_5$ and consider the function
$\Lambda:[0,1]\to{\mathcal C}^{\infty}\left([1,\infty),\R\right)$ that assigns $\mu\in[0,1]$ to 
$\Lambda(\mu)=\mu\tau_4+(1-\mu)\rho.$
Let us also take $q_1=q_2$ so that 
\begin{equation}\label{astuce}
\Lambda(\lambda_4)'(q_1)=\tau'(q_1)
\end{equation}
(the parameter $\lambda_4$ depends on $q_2$). It is sufficient to show that $\Lambda(\lambda_4)\le\tau$ : to obtain that $\Lambda(\lambda_4)\le\tau$ outside $[q_1,q_2]$, for  $q_1<q_2$,  one can use a simple continuity  argument on the graph of $\Lambda(\lambda_4)$, seen as a function of $q_2$.

By (\ref{astuce}) we obtain $\displaystyle \frac{(\tau-\rho)'(q_1)}{(\tau_4-\rho)'(q_1)}=1>\lambda_4$. The function  
$\displaystyle\frac{\tau-\rho}{\tau_4-\rho}$ being increasing in a neighborhood of $q_1$  (as we will show below) we get that for $q>q_1$
$(\tau-\rho)(q)>\lambda_4(\tau_4-\rho)(q)$ which implies  $\Lambda(\lambda_4)(q)<\tau(q)$ for $q\not=q_1$.

To finish the proof we need to show that $\displaystyle\frac{\tau-\rho}{\tau_4-\rho}$ is increasing. 
Put $s_1=\frac{\lambda_3}{\lambda_3+\lambda_5}-\lambda_1$, 
$k_5=\frac{\lambda_5}{\lambda_3+\lambda_5}=1-k_1$ and $k_1=\frac{\lambda_3}{\lambda_3+\lambda_5}$, all positive. We can write:
\begin{eqnarray*}
\frac{\rho-\tau}{\rho-\tau_4}&=&\frac{s_1(\tau_1-\tau_2)+k_5(\tau_5-\tau_2)}{k_1(\tau_1-\tau_4)+k_5(\tau_5-\tau_4)}=\\
&=&\frac{1}{1-\lambda_1}\frac{\displaystyle s_1-\frac{\tau_2-\tau_5}{\tau_1-\tau_5}}{\displaystyle k_1-\frac{\tau_4-\tau_5}{\tau_1-\tau_5}}=\frac{k_1}{s_1(1-\lambda_1)}\frac{1-f}{1-g}
\end{eqnarray*}
where $f,g$  are both positive increasing functions by lemma \ref{subsidiary}. Moreover $f/g$ is increasing which implies $f'g-g'f>0$ and $g'(q_1)=f'(q_1)$ hence $(f-g)'(g+f)+g'-f'<0$. The result follows.
\end{Proof}

The proof of theorem \ref{cocorico} in the case $n>2$ is now easy to derive : suppose $\tau=\sum_{k=1}^{n}\lambda_k\tau(p_k,.)$ and let $\tau(p_1,.)$ and $\tau(p_2,.)$ be the first two functions of the convex combination. Be the previous two lemmas there exist a convex combination $\hat\tau$ of three $\tau(p_i,.)$ functions such that 
\begin{enumerate}
\item $\frac{1}{\lambda_1+\lambda_2}\left(\lambda_1\tau_1(q_i)+\lambda_2\tau_2(q_i)\right)=\hat\tau(q_i)$ , for $i=1,2$  
\item $\frac{1}{\lambda_1+\lambda_2}\left(\lambda_1\tau_1'(q_1)+\lambda_2'\tau_2(q_1)\right)<\hat\tau'(q_1)$ \; , \; $\frac{1}{\lambda_1+\lambda_2}\left(\lambda_1\tau_1'(q_2)+\lambda_2'\tau_2(q_2)\right)>\hat\tau'(q_2)$ 
\item $\frac{1}{\lambda_1+\lambda_2}\left(\lambda_1\tau_1+\lambda_2\tau_2\right)<\hat\tau$ on $(q_1,q_2)$ and $\frac{1}{\lambda_1+\lambda_2}\left(\lambda_1\tau_1+\lambda_2\tau_2\right)>\hat\tau$  on $(1,\infty)\setminus[q_1,q_2]$. 
\end{enumerate}
The function $\tilde\tau=(\lambda_1+\lambda_2)\hat\tau+\sum_{k=3}^{n}\lambda_k\tau(p_k,.)$ satisfies then the conclusion of theorem \ref{cocorico}.

We can now prove theorem \ref{3} :

{\em There exists an inhomogeneous Bernoulli product $\mu$ such that the spectrum $\tau$ of $\mu$ is not derivable on a dense subset of \ $[1,\infty)$. }

The strategy of the demonstration of this theorem is the following : we first find inhomogeneous Bernoulli products that are not derivable at a finite number of predefined points and we construct the measure $\mu$ using Cantor's diagonal argument.
\begin{Lemma}
For any $p_1,...,p_n$ and any convex combination $\tau$ of $\tau(p_1,.),...,\tau(p_n,.)$ there exist an inhomogeneous Bernoulli measure $\mu$ whose multifractal spectrum equals $\tau$.
\end{Lemma}
The proof of this lemma is not difficult and left to the reader. Let un now prove theorem \ref{3}.
\begin{Proof}{\bf of Theorem \ref{3}.}
Fix $(q_n)_n$ a sequence of real numbers, dense in $[1,\infty)$ and nested in the sense that $q_{2n+1}<q_{2n+2}$ and $\{q_1,...,q_{2n}\}\cap[q_{2n+1},q_{2n+2}]=\emptyset$ for all $n\ge 0$.  Let $p_1,p_2\in (0,1)$ and $\tau_1=\frac12\tau(p_1,.)+\frac12\tau(p_2,.)$. By the previous lemma we can construct a Bernoulli product $\mu_1$ of spectrum $\tau_1$. Theorem \ref{cocorico} implikes then the existence of a convex combination $\tau_2$ of $\tau(p_i,.)$'s functions, such that 
\begin{enumerate}
\item $\tau_1(q_i)=\tau_2(q_i)$ , for $i=1,2$ , $\tau_1'(q_1)<\tau_2'(q_1)$ \; , \; $\tau_1'(q_2)>\tau_2'(q_2)$  
\item $\tau_2>\tau_1$ on $(q_1,q_2)$ and  $\tau_2<\tau_1$  on $(1,\infty)\setminus[q_1,q_2]$. 
\end{enumerate}
We can therefore define a measure $\mu_2$ of spectrum $\tau_2$; Using $\mu_1$ and $\mu_2$ we can construct a measure $\nu_2$ of spectrum $\max\{\tau_1,\tau_2\}$: Let $\mu_i$ me the inhomogeneous Bernoulli measure of spectrum $\tau_i$, $i=1,2$. Take $(\ell_k)_k$ a sequence of integers such that $\displaystyle\frac{\ell_{k+1}}{\sum_1^k\ell_i}\to\infty$. On dyadique intervalles of length between $2^{-\ell_{2k}}$ and $2^{-\ell_{2k+1}}$ apply the weight distribution of $\mu_1$ and on dyadique intervalles of length between $2^{-\ell_{2k+1}}$ and $2^{-\ell_{2k+2}}$ apply the weight distribution of $\mu_2$, where $k\in \N$. It is easy to verify that the resulting inhomogeneous measure $\nu_2$ has spectrum $\max\{\tau_1,\tau_2\}$. Remark that this spectrum equals $\tau_2$ on $[q_1,q_2]$ and $\tau_1$ elsewhere on $[1,\infty)$.

We proceed by induction. Suppose the measures $\nu_1=\mu_1$, $\mu_2$, $\nu_2$,..., $\mu_n$, $\nu_n$ defined and denote $\tau_i$ the spectrum of the measure $\mu_i$, $i\in\{1...,n\}$. We assume that that on every interval $[q_{2i+1},q_{2i+2}]$ , where $i\le n$, the spectrum of $\nu_n$ equals  $\max\{\tau_1,...\tau_n\}$ and is realized by $\tau_i$. Let us construct $\mu_{n+1}$ and $\nu_{n+1}$. Consider $\tau_j$ the function that equals the $\max\{\tau_1,...\tau_n\}$ on  $[q_{2(n+1)+1},q_{2(n+1)+2}]$. By theorem \ref{cocorico} we can find a function $\tau_{n+1}$ satisfying : 
\begin{enumerate}
\item $\tau_{n+1}(q_{2(n+1)+i})=\tau_j(q_{2(n+1)+i})$ , for $i=1,2$ ,\\
 $\tau_{n+1}'(q_{2(n+1)+1})>\tau_j'(q_{2(n+1)+1})$ \; , \; $\tau_{n+1}'(q_{2(n+1)+2})<\tau_j'(q_{2(n+1)+2})$  
\item { $\tau_{n+1}>\tau_j \mbox{ on }(q_{2(n+1)+1},q_{2(n+1)+2})$ and} 
{ $\tau_{n+1}<\tau_j\mbox{  on }(1,\infty)\setminus[q_{2(n+1)+1},q_{2(n+1)+2}].$}
\end{enumerate}
Let $\mu_{n+1}$ be the inhomogeneous Bernoulli measure of spectrum $\tau_{n+1}$.  To define the measure $\nu_{n+1}$ we use the previous procedure convenably adapted:
Take $(\ell_k)_k$ a sequence of integers such that $\displaystyle\frac{\ell_{k+1}}{\sum_1^k\ell_i}\to\infty$. On dyadique intervalles of length between $2^{-\ell_{(n+1)k+i}}$ and $2^{-\ell_{(n+1)k+i+1}}$ apply the weight distribution of $\mu_i$, where $i=1,...,n+1$, $k\in \N$. It is easy to verify that the resulting inhomogeneous measure $\nu_{n+1}$ has spectrum $\max\{\tau_1,...\tau_{n+1}\}$ on $(1,\infty)$. Remark that this spectrum equals $\tau_{n+1}$ on $[q_{2(n+1)+1},q_{2(n+1)+2}]$ and $\max\{\tau_1,...\tau_n\}$ elsewhere on $[1,\infty)$.

To end the proof we use Cantor's diagonal argument: take $(\ell_k)_k$ comme ci-dessus et  define the measure $\nu$ by applying the weight distribution of $\nu_k$ on dyadique intervalles of length between $2^{-\ell_{k}}$ and $2^{-\ell_{k+1}}$. The spectrum of the measure $\nu$ equals then $\tau=\sup_{n\in\N}{\tau_n}$. By  construction the function $\tau$ is not derivable at the points $(q_k)_k$ and the proof of theorem \ref{3} is complete. 
\end{Proof}
\bibliographystyle{alpha}
\bibliography{biblio}
\end{document}